\documentclass[12pt]{article}
\usepackage{a4wide}
\usepackage{epsfig}
\usepackage{amsmath,amsfonts,amssymb,amsthm}
\usepackage{eucal}
\usepackage{graphics,graphicx}
\usepackage{float}
\usepackage{subfig}
\usepackage{epstopdf}
\usepackage{authblk}
\newtheorem{thm}{Theorem}[section]

\newtheorem{lem}{Lemma}[section]
\newtheorem{rem}{Remark}[section]
\newcommand{\C}{\mathbb{C}}

\begin{document}

\title{{Different types of wandering domains in the family $ \lambda+z+\tan z$}}

\author{Subhasis Ghora\footnote{sg36@iitbbs.ac.in(Corresponding author)}} 
	
%\author[1]{Tarakanta Nayak\footnote{tnayak@iitbbs.ac.in} }
%\affil[1]{\textit{School of Basic Sciences \hspace{10cm}Indian Institute of Technology Bhubaneswar, India}  }
\affil{\textit{Department of Mathematics \hspace{10cm}C. V. Raman Global University,  Bhubaneswar, India}  }
\date{}

%\author{Subhasis Ghora\\ 
%Department of Mathematics, C. V. Raman Global University}

\date{}
%%%%%%%%%%%%%%%%%%%%%%%%%%%%%%%%%%%%%%%%%%%%%%%%%%%%%%%%%%%%%%
%\setcounter{page}{0}
%\begin{center}
%
%\begin{tabular}{|p{5cm}c| p{9cm}|} \hline
%Title of the Paper: & & Iteration of $\lambda+\tan z$\\ \hline Name
%of the Authors: & & Tarakanta Nayak
%\\ \hline
%Affiliations: & & \parbox[t]{8cm}{

%Tarakanta Nayak \\
%Ph.~D. Student \\
%Department of Mathematics \\
% Indian Institute of Technology Guwahati\\
% Guwahati - 781039, India\\}
%  \\ \hline
%  Proposed Running Title: & &  Iteration of $\lambda+\tan z$ \\ \hline
%   Mathematics Subject Classification 2000~ Primary & & 37F45, 37F50\\
%\hline
% Key  words &
% &  Locally connected, Jordan curve, Julia set, Transcendental meromorphic functions \\
%\hline
%Name and Address of the Author to contact: &  & \parbox[t]{8cm}{M.~Guru Prem Prasad \\
%Assistant Professor \\
%Department of Mathematics \\
% Indian Institute of Technology Guwahati\\
% Guwahati - 781039, India\\
%Phone: +91-361-2582608 \\
%Fax: +91-361-2690762 \\
%Email: mgpp@iitg.ernet.in } \\ \hline
%\end{tabular}
%\end{center}
%\clearpage

%%%%%%%%%%%%%%%%%%%%%%%%%%%

\maketitle
\begin{abstract}
Dynamics of an one-parameter family of functions $f_\lambda(z)=\lambda + z+\tan z, z \in \mathbb{C}$ and $\lambda \in \mathbb{C}$ with an unbounded set of singular values is investigated in this article. 	For $|2+\lambda^2|<1$, $\lambda=i$, $2+\lambda^2=e^{2\pi i \alpha}$ for some rational number $\alpha$ and for some bounded type irrational number $\alpha$, the dynamics of $f_{\lambda+m\pi}$ is determined for $m \in \mathbb{Z}\setminus\{0\}$. For such values of $\lambda$, the existence of $m$ many wandering domains of $f_{\lambda+m\pi}$ with disjoint grand orbits in the lower half-plane are asserted along with a completely invariant Baker domain containing the upper half-plane. Further, each of such wandering domains is found to be simply connected, unbounded, and escaping. Different types of the internal behavior of $\{f^n_{\lambda+m\pi}\}$ on such a wandering domain $W$ are highlighted for different values of $\lambda$. More precisely, for $\mid2+\lambda^2\mid<1$, it is manifested that the forward orbit of any point $z\in W$ stays away from the boundaries of $W_n$s. For $\lambda=i$, it is proved that  $\liminf_{n\rightarrow \infty}dist(f^n_{i+m\pi}(z),\partial W_n)=0$ for all $z\in W$. Further, $\Im(f^n_{i+m\pi}(z))\rightarrow -\infty$ as $n \rightarrow \infty$. For $2+\lambda^2=e^{2\pi i\alpha}$ for some rational number $\alpha$, $\liminf_{n\rightarrow \infty}dist(f^n_{\lambda+m\pi}(z),\partial W_n)=0$  is established for all $z\in W$. But, $\Im(f^n_{\lambda+m\pi}(z))$ tends to a finite point for all $z\in W$ whenever $n \rightarrow \infty$. For $2+\lambda^2=e^{2\pi i\alpha}$, $\liminf_{n\rightarrow \infty}dist(f^n_{\lambda+m\pi}(z),\partial W_n)>0$ for all $z\in W$ and $dist(f^n_{\lambda+m\pi}(z),f^n_{\lambda+m\pi}(z')=dist(z,z')$ is authenticated for all $z,z'\in W$ and for some bounded type irrational number $\alpha$.

	2010 Mathematics Subject Classification 37F50  30D05  30D030 
\end{abstract}

\section{Introduction}
A transcendental meromorphic map $f:\mathbb{C} \rightarrow \widehat{\mathbb{C}}$ is referred to as general meromorphic if it has either exactly one pole that is not an omitted value or at least two poles. For the sake of simplicity, we refer to such a general meromorphic map in this article as a meromorphic map. Depending upon the behavior of the sequence of functions $\{f^n\}$ in the complex plane, the extended complex plane can be fissioned into two sets namely, the Fatou set and the Julia set. The Fatou set consists of all the points in the complex plane for which $\{f^n\}$ is defined and normal in a neighborhood of that point. The Julia set is the complement of the Fatou set in $\widehat{\mathbb{C}}$.
\par 
 A Fatou component is a maximal connected subset of the Fatou set. A Fatou component $U$ of a function $f$ is called $p$-periodic if $p$ is the smallest natural number such that $f^p(U)\subset U$. The set $\{U, U_1, U_2,...,U_{p-1}\}$ is the cycle containing $U$ where $U_n=f^n(U)$. The Fatou component, $U$ is called invariant if $p=1$. If $f(U)\subset U$ and $f^{-1}(U)\subset U$ then $U$ is called completely invariant. A periodic Fatou component $U$ can be categorized as an attractive domain, a parabolic domain, a Baker domain, a Herman ring, or a Siegel disc, based on the internal behavior of $\{f^n\}$ on $U$. For each pair of distinct natural numbers $m,n$, if $f^m(U)\neq f^n(U)$ then $U$ is called wandering.
 
 \par
  
 The image of a critical point, i.e., $f(z_{0})\ \mbox{where}\ f'(z_{0})=0$, is known as a critical value of $f$. A point $a\in \widehat{\mathbb{C}}$ is an asymptotic value of $f$ if there exists a curve $\gamma:[0,\infty)\rightarrow \mathbb{C}$ with $\lim_{t\rightarrow \infty} \gamma(t)=\infty$ such that $ \lim_{t\rightarrow \infty} f(\gamma(t))=a$. If $D_r(a)$ is a disk (with respect to the spherical metric) centered at an asymptotic value $a \in \widehat{\C}$  of $f$ of radius $r$ then a component $V_r$ of $f^{-1}(D_r(a))$ can be plumped for in such a way that $V_{r_1}\subset V_{r_2}$ for $0<r_1<r_2$. In this case, $\bigcap_{r>0} V_r= \emptyset$ and the choice $r \mapsto V_r$ defines a transcendental singularity. We say that a singularity $V$ lies over $a$. A singularity is called logarithmic if  $f: V_r \to D_{r}(a) \setminus \{a\}$ is a universal covering. In this case, $V_r$ is simply connected. It is called non-logarithmic if $V_r$ is not simply connected.  The set of singular values is defined as the closure of all critical values and asymptotic values of $f$. It is denoted by $S_f$.

 \par
 Studying the dynamics of particular functions has been incredibly helpful in forecasting outcomes for a larger class of functions containing them as well as frequently providing hints for the proofs of general theorems. The study of the dynamics of meromorphic maps is largely developed around the functions with the bounded set of singular values \cite{keen-kot, sajid3, sajid4}. To our knowledge, the dynamics of the family of meromorphic maps with an unbounded set of singular values is only investigated in~\cite{nayakmgp}. In light of the substantial body of literature on the dynamics of functions with a bounded set of singular values and the dynamics of functions with a logarithmic type of singularity, it will not be unwise to conclude that the study of the dynamics of functions with an unbounded set of singular values is essentially a new arena of research and mostly remains unexplored. For such functions, if there is at least one singular value over which there is a non-logarithmic singularity then the family becomes more prominent and worth studying the dynamics. Imbued by this, we investigate the dynamics of a family of functions namely,
$$f_\lambda(z)=\lambda+z+\tan z$$ for $\lambda\in \mathbb{C}$.  Note that $f_{-\lambda}(z)=-\lambda+z+\tan z=-f_\lambda(-z)$ concluding that $f_\lambda$ and $f_{-\lambda} $ are conformally conjugate and hence, the dynamical behavior of $f_\lambda$ and $f_{-\lambda}$ are analogous. Thus, it is enough to study the dynamics of $f_\lambda$ for $\Im(\lambda)>0$. It is proved in \cite{tanz1} that $f_\lambda$ has an unbounded set of singular values and there is a non-logarithmic singularity lying over $\infty$. They also have proved that for $\Im(\lambda)>0$, $f_\lambda$ has a completely invariant Baker domain containing the upper-half plane. It is referred to as the {\it{primary Fatou component}}. For some values of $\lambda$ for which $(\overline{\cup_{s\in S_{f_\lambda}}\cup_{n\geq0}f_\lambda^n(s)})\cap \mathbb{C}$ is contained in the Fatou set, the dynamics of $f_\lambda$ is investigated in  \cite{tanz1}. This, however, is not always the case in this article. We have investigated the dynamics of $f_{\lambda+m\pi}$ for all $\lambda$ such that $|2+\lambda^2|<1$ or $\lambda=i$ or $2+\lambda^2=e^{2\pi i \alpha}$ for some $\alpha$ rational or $2+\lambda^2=e^{2\pi i \alpha}$ for some bounded type irrational number $\alpha$. Based on the internal behavior of the sequence of functions on a simply connected wandering domain of an entire function, \cite{classi} provides a classification of such wandering domains. One of such classifications is based on the behavior of the orbit of a point on such wandering domains to their boundaries. Several examples of entire functions are given to visualize different types of behavior. This class of meromorphic maps likewise exhibits these internal behaviors, as can be seen in this article. This family is so dynamically rich that the presence of all kinds of such behaviors is perceived in this family.

\par
For any non-zero integer $m$ and $\lambda$ such that $\mid2+\lambda^2\mid<1$, the dynamics of $f_{\lambda+m\pi}$ is completely determined and is demonstrated in the following theorem. The existence of escaping the wandering domain is proved for such values of $\lambda+m\pi$. It is shown that the wandering domains are unbounded and simply connected. The precise statement is the following.
\begin{thm}\label{wanderingattracting}
For $\mid2+\lambda^2\mid<1$, in addition to the primary Fatou component, $f_{\lambda+m\pi}$   has $m$ many  wandering domains with distinct grand orbits. If $W$ is such a wandering domain then  it has the following properties.

\begin{enumerate}
\item Each $W$ is escaping.

\item There is a two-sided sequence of unbounded wandering domains 
$\{W_n \}_{n \in \mathbb{Z}}$ in the grand orbit of $W$ such that  $f_{\lambda+m\pi}: W_{n} \to W_{n+1}$ is a proper map with degree $2$. 
\item If  $W'$ is a wandering domain in the grand orbit of $W$ and  different from all $W_n$s then  $f_{\lambda+m\pi}$ is one-one on $  W'  $.
\item $\liminf_{n\rightarrow \infty}dist(f_{\lambda+m\pi}^n(z),\partial W_n)>0$ for all $z\in W$, i.e., the forward orbit of $z$ stays away from the boundary.
\end{enumerate} 
The Fatou set of $f_{\lambda+m\pi}$ is the union of the primary Fatou component and these $m$ many grand orbits of wandering domains. 
\end{thm}
 For $\lambda=i+m\pi$, the dynamics of $f_{i+m\pi}$ is elucidated in the following theorem. The internal behavior of the wandering domains of $f_{i+m\pi}$ is found to be different from that of the wandering domains described in Theorem \ref{wanderingattracting}.
 \begin{thm}\label{wanderingbaker}
 	For $\lambda=i$, in addition to the primary Fatou component, $f_{i+m\pi}$   has $m$ many wandering domains with distinct grand orbits. If $W$ is such a wandering domain then it has the following properties.
 	
 	\begin{enumerate}
 		\item Each $W$ is escaping.
 		
 		\item There is a two-sided sequence of unbounded wandering domains
 		$\{W_n \}_{n \in \mathbb{Z}}$ in the grand orbit of $W$ such that  $f_{i+m\pi}: W_{n} \to W_{n+1}$ is a proper map with degree $2$.
 		\item If  $W'$ is a wandering domain in the grand orbit of $W$ and different from all $W_n$s then  $f_{i+m\pi}$ is one-one on $  W'  $.
 		\item $\liminf_{n\rightarrow \infty}dist(f^n_{i+m\pi}(z),\partial W_n)=0$ for all $z\in W$. Further, $\Im(f^n_{i+m\pi}(z))\rightarrow -\infty$ for all $z\in W$ whenever $n \rightarrow \infty$.
 	\end{enumerate}
 	The Fatou set is the union of the primary Fatou component and these $m$ many grand orbits of wandering domains.
 \end{thm}

 All the fixed points of $f_{\lambda}$ (if exists) are indifferent if  $2+\lambda^2=e^{2\pi i \alpha}$. These fixed points are parabolic if $\alpha$ is rational except for $\lambda=i$. The next theorem explores the dynamics of $f_{\lambda+m\pi}$ for all such values of $\lambda$.
 \begin{thm}\label{wanderingparabolic}
 	For $2+\lambda^2=e^{2\pi i\alpha}$ for some rational number $\alpha$ and $m(\neq 0)\in \mathbb{Z}$, in addition to the primary Fatou component, $f_{\lambda+m\pi}$   has $m$ many  wandering domains with distinct grand orbits. If $W$ is such a wandering domain then it has the following properties.
 	
 	\begin{enumerate}
 		\item Each $W$ is escaping.
 		
 		\item There is a two-sided sequence of unbounded wandering domains
 		$\{W_n \}_{n \in \mathbb{Z}}$ in the grand orbit of $W$ such that  $f_{\lambda+m\pi}: W_{n} \to W_{n+1}$ is a proper map with degree $2$.
 		\item If  $W'$ is a wandering domain in the grand orbit of $W$ and different from all $W_n$s then  $f_{\lambda+m\pi}$ is one-one on $  W'  $.
 		\item $\liminf_{n\rightarrow \infty}dist(f^n_{\lambda+m\pi}(z),\partial W_n)=0$ for all $z\in W$, i.e., the forward orbit of $z$ converges to the boundary. But $\Im(f^n_{\lambda+m\pi}(z))$ tends to a finite point for all $z\in W$ whenever $n \rightarrow \infty$
 	\end{enumerate}
 \end{thm}

An irrational number is called bounded type if its continued fraction coefficients are bounded. It is known that if the multiplier of a fixed point $z_0$ of a function $f$ is $e^{2\pi i \alpha}$ where $\alpha$ is a bounded type irrational number then $z_0$ is a Siegel point \cite{milnor}. In other words, $f$ has a Siegel disc centered at $z_0$. Note that for $2+\lambda^2=e^{2\pi i\alpha}$ for some bounded type irrational $\alpha$, $f_\lambda$ has infinitely many Siegel points and hence, infinitely many Siegel discs. We describe the dynamics of $f_{\lambda+m\pi}$ where $2+\lambda^2=e^{2\pi i\alpha}$ for some bounded type irrational number $\alpha$.

\begin{thm}\label{wanderingsiegel}
For $2+\lambda^2=e^{2\pi i\alpha}$ for some bounded type irrational $\alpha$ and $m(\neq 0)\in \mathbb{Z}$, in addition to the primary Fatou component, $f_{\lambda+m\pi}$   has $m$ many  wandering domains with distinct grand orbits. If $W$ is such a wandering domain then  it has the following properties.

\begin{enumerate}
\item Each $W$ is escaping. 
\item For any wandering domain $W$,  $f_{\lambda+m\pi}:W\rightarrow W_1$ is one-one.
\item $\liminf_{n\rightarrow \infty}dist(f^n_{\lambda+m\pi}(z),\partial W_n)>0$ for all $z\in W$ and $dist(f^n_{\lambda+m\pi}(z),f^n_{\lambda+m\pi}(z')=dist(z,z')$ for all $z,z'\in W$.
\end{enumerate} 
\end{thm}

\par

 For any set $A\subset \C$, the closure of $A$ is denoted by $\overline{A}$ and the boundary of $A$ is denoted by $\partial A$. Let $D_r(a)$ denote the disc centered at $a$ of radius $r$. The Euclidean distance of any two points $z,~w$ (or two sets $A$, $B$) is  denoted by $dist(z,w)$ (or $dist(A,B)$ respectively). The imaginary and the real part of a complex number $z$ is denoted by $\Im(z)$ and $\Re(z)$ respectively.

%For a curve $\gamma$, $l(\gamma)$ denotes the length of $\gamma$. The set of all parabolic fixed points of $f$ is denoted by $\mathcal{P}_{f}$.

\section{Preliminaries}

The following lemma demonstrated in \cite{bolsch} is to be used frequently. If the pre-image of any compact subset of $U$ is compact in $V$ for a continuous map $f: V \rightarrow U$ between two open connected subsets of $\mathbb{C}$, then $f$ is considered proper. Furthermore, if $f$ is analytic, then there is a $d$ such that every element of $U$ has $d$ pre-images when multiplicity is taken into account. In this case, the local degree of $f$ at a point $z$ equals the multiplicity of that point. This value $d$ is referred to as the degree of $f: V \to U$.

\begin{lem}\label{RH}{\bf{(Riemann-Hurwitz formula)}}

Let $f:\mathbb{C}\rightarrow \widehat{\mathbb{C}}$ be a transcendental meromorphic function. If $V$ is a  component of the pre-image of  an open connected set $U$ and $f : V \rightarrow U$ is a proper map of degree $d$, then $c(V)-2=d(c(U)-2)+n$, where $n$ is the number of critical points of $f$ in $V$ counting multiplicity and $n \leq 2d-2$. Here, the multiplicity of a critical point is one less than the local degree of $f$ at the critical point.  
\end{lem}

Let $f_\lambda(z)=\lambda+z+\tan z$ for  $\lambda \in \mathbb{C}$. Several dynamical aspects of $f_\lambda$ are ascertained in \cite{tanz1} for different values of $\lambda$. As a prelusory, we describe a few of them in Lemma \ref{criticalpoint}, Lemma \ref{fixedpoint}, Lemma \ref{cifcupperhalf} and Lemma \ref{attractingdomaincomplete}. The set of all singular values of $f_\lambda$ is also found. It is seen that $f_\lambda$ has an unbounded set of critical values contained in two horizontal lines, one in the upper half-plane and the other in the lower half-plane. The point at $\infty$ is the only asymptotic value. A detailed description of the singular values of $f_\lambda$ is described below.

\begin{lem}\label{criticalpoint}{\cite{tanz1}}
	\begin{enumerate}
\item The set of all critical points of $f_\lambda$ is  $\{\frac{\pi}{2}+n \pi  \pm i \sinh^{-1}1 : n \in \mathbb{Z}\}$ and  $ \lambda+\frac{\pi}{2}+n\pi\pm i(\sinh^{-1}1+\sqrt{2})$ are the critical values of $f_\lambda$.
\item  $f_\lambda$ has only one asymptotic value. The asymptotic value is the point at infinity. Further, the singularity lying over the point at infinity is non-logarithmic type.  

	\end{enumerate}
\end{lem}

For $\lambda=i$, $f_\lambda$ does not have any fixed point. For all other $\lambda$s with $\Im(\lambda)>0$, there are infinitely many fixed points of $f_\lambda$. It is important to note that for a given $\lambda$, the multiplier of all the fixed points of $f_\lambda$ is equal.

\begin{lem}{\cite{tanz1}}\label{fixedpoint}
	\begin{enumerate}
\item  All the fixed points of $f_\lambda$ are in the lower half-plane whenever $\Im(\lambda)>0$. A point $z$ is a fixed point of $f_\lambda$ if and  only if $z+n \pi$ is so for all $n \in \mathbb{Z}$.
		\item The multiplier of each fixed point is $2+\lambda^2$. Hence, all the fixed points of $f_\lambda$ are attracting,  repelling, or indifferent  together.

	\end{enumerate}
\end{lem}

The existence of a completely invariant Baker domain of $f_\lambda$ is also established in \cite{tanz1} for all $\lambda$. This completely invariant Baker domain of $f_\lambda,$ is referred to as the \textit{ primary Fatou component}. This is described in the following lemma.

\begin{lem}{\cite{tanz1}}\label{cifcupperhalf}
	For  $\Im(\lambda) > 0,$ there is a   completely invariant Baker domain of $f_\lambda$  containing the upper half-plane.  
\end{lem}

For all $\lambda$ with $\mid2+\lambda^2\mid<1$, it is seen that $f_\lambda$ has infinitely many attracting fixed points. The dynamics of $f_\lambda$ for such values of $\lambda$ are described in the following lemma.
\begin{lem}{\cite{tanz1}}\label{attractingdomaincomplete}
	Let $\mid2+\lambda^2\mid<1$. Then, along with the primary Fatou component, there are infinitely many invariant attracting domains of $f_{\lambda}$ and each such attracting domain $U$ is unbounded in such a way that $\{\Im(z): z \in U\}$ is unbounded but for every $z_0\in U$,  $\{\Re(z): z \in U ~\mbox{and}~\Im(z)>\Im(z_0)\}$ is bounded and $f_\lambda:U \rightarrow U$ is a proper map of degree $2$. Further, $f_{\lambda}$ does not have any other periodic Fatou component or any wandering domain.
\end{lem}

In \cite{access, fagella2019}, the dynamics of $f_\lambda$ for $\lambda=i$ is portrayed. The existence of infinitely many invariant Baker domains is found. It is noted below.

\begin{lem}\label{bakerdomaincomplete}
	For $\lambda=i$, along with the primary Fatou component, there are infinitely many invariant Baker domains of $f_{\lambda}$ and each such Baker domain $U$ is unbounded in such a way that $\{\Im(z): z \in U\}$ is unbounded but,  $\{\Re(z): z \in U \}$ is bounded and $f_\lambda:U \rightarrow U$ is a proper map of degree $2$. Further, $f_{\lambda}$ does not have any other periodic Fatou component or any wandering domain.
\end{lem}

Now we are going to describe the dynamics of $f_\lambda$ for all $\lambda \in \mathbb{C}$ such that $2+\lambda^2=e^{2\pi i \alpha}$ for some rational number $\alpha$. In this case, $f_\lambda$ has infinitely many parabolic fixed points in the lower half-plane. Corresponding to each such fixed point,  the Fatou set of $f_\lambda$ has infinitely many parabolic domains, described in the following lemma.

\begin{lem}\label{parabolic}
	For $2+\lambda^2=e^{2\pi i \alpha}$, where $\alpha$ is rational, in addition to the primary Fatou component, there are infinitely many unbounded parabolic domains in the lower half plane. Moreover, there is no other periodic Fatou component except possibly the Baker domain.
\end{lem}

\begin{proof}
	
Each fixed point of $f_\lambda$ satisfies the equation $\tan z=-\lambda$. Since $\Im (\lambda)>0$ and $\lambda\neq i$, there are infinitely many fixed points in the lower half-plane. Note that if $z_0$ is a fixed point of $f_\lambda$ then for each $m\in \mathbb{Z}$, $z_0+m\pi$ is also a fixed point of $f_\lambda$. As the multiplier of each of the fixed points is $2+\lambda^2=e^{2\pi i \alpha}$ for some rational number $\alpha$, each of the fixed points is parabolic. There are infinitely many invariant parabolic domains corresponding to such fixed points. Each such domain is simply connected. Since $f_\lambda$ does not have any weakly repelling fixed points, it follows from Theorem B of \cite{access} that such parabolic domains are unbounded. If $c$ is a critical point in the lower half-plane such that $f_\lambda^n(c)\rightarrow z_0$ then $f_\lambda^n(c+m\pi)\rightarrow z_0+m\pi$. Thus all the critical points in the lower half-plane are contained in such parabolic domains. Also, all the critical points in the upper half-plane are contained in the primary Fatou component. This gives that $f_\lambda$ does not have any other attracting domains,  parabolic domains, or Siegel discs. The existence of the primary Fatou component rules out the existence of the Herman ring. This completes the proof.

	\begin{figure}[H]
		\centering
		\subfloat[Parabolic domain:  $ \alpha= \frac{1}{8}$]
		{\includegraphics[width=2.5in,height=2.5in]{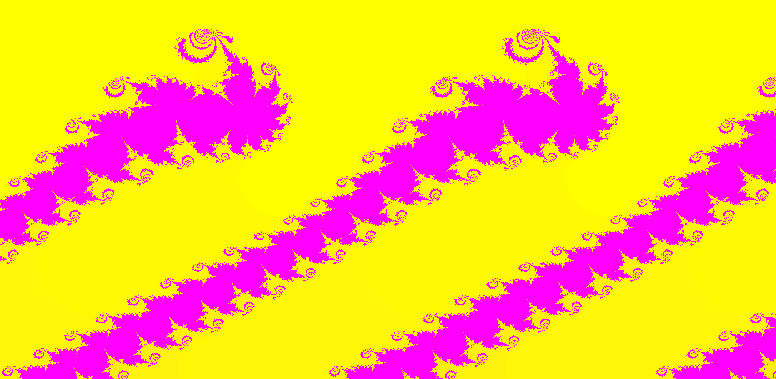}}
		\vspace{.2cm}
		\subfloat[Magnification of one of the parabolic domain]	{\includegraphics[width=2.5in,height=2.5in]{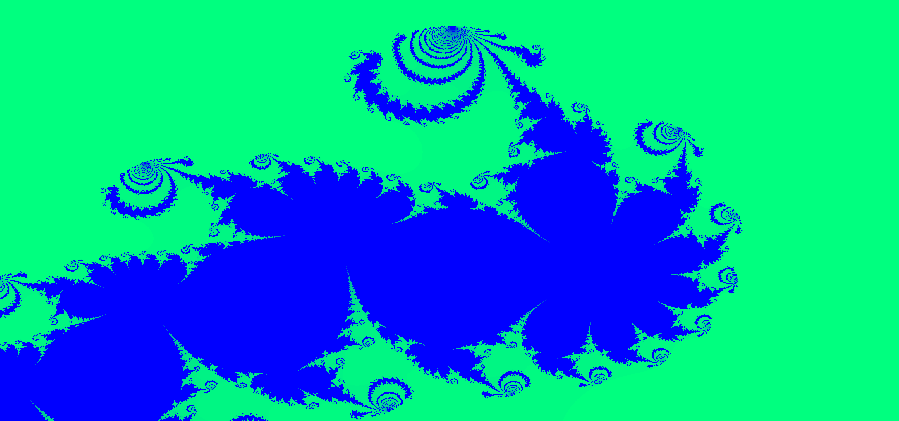}}
		\caption{ Parabolic domains for $\lambda=\frac{\sqrt{2(-4+\sqrt{2}+i\sqrt{2})}}{2}$} or  $ \alpha= \frac{1}{8}$.	 
		\label{Fatouset-3}
	\end{figure}
	
\end{proof}

\begin{rem}\label{ Siegel}
	Note that each fixed point of $f_\lambda$ is a  Siegel point whenever  $2+\lambda^2=e^{2\pi i\alpha}$ for some bounded type irrational $\alpha$. In that case, all the critical points of $f_\lambda$ are in the Julia set. Further, in addition to the primary Fatou component, $f_\lambda$ has infinitely many  Siegel discs in the lower half plane. Moreover, there is no other attracting domain, parabolic domain, and Herman ring of $f_\lambda$.
\end{rem}

\begin{rem}\label{translationofad}(Enumeration of the Fatou components)
	
	\begin{enumerate}
		\item For a Fatou component $V$ of $f_\lambda$ and $k\in \mathbb{Z}$, let $V+k\pi=\{z+k\pi:z\in V\}$. As $f_\lambda^n(z+k\pi)=k\pi +f_\lambda^n(z)$, then $z$ is in the Fatou set of $f_\lambda$ if and only if $z+k\pi$ is also in the Fatou set of $f_\lambda$. To put it another way, the Fatou set is $\pi$-invariant implying that for each $k\in \mathbb{Z}$, $V+k\pi$ is a Fatou component of $f_\lambda$. 
		
		\item  For $k_1\neq k_2$, if $V+k_1\pi \cap V+k_2\pi \neq \emptyset$ then $\cap_{k\in \mathbb{Z}}V+k\pi\neq \emptyset$. In fact, $\cap_{k\in \mathbb{Z}}V+k\pi =V$. As $V+k_1\pi \cap V+k_2\pi \neq \emptyset$, for any $z\in V$, there exists a curve $\gamma$ contained in the Fatou set connecting $z+k_1\pi$ to $z+k_2\pi$. Since the Fatou set is $\pi$-invariant and the distance between $z+k_1\pi$ and $z+k_2\pi$ is greater than or equals to $\pi$, for each $k\in \mathbb{Z}$, $\gamma+k\pi$ is in the Fatou set. By the construction of $V+k\pi$, $\cup_{k\in \mathbb{Z}}\gamma+k\pi$ intersects each $V+k\pi$ giving that  $\cap_{k\in \mathbb{Z}}V+k\pi\neq \emptyset$. As two Fatou components are either identical or disjoint,  $\cap_{k\in \mathbb{Z}}V+k\pi =V$. Thus for a given Fatou component $V$ of $f_\lambda$, either each of $V+k\pi$ is disjoint or all of them coincide with $V$.
	\item 
 Let $V$ be a $p$-periodic Fatou component of $f_\lambda$ and $\{V, V^1, V^2,..., V^{p-1}\}$ be the cycle containing $V$, where $V^i=f_\lambda^i(V)$. We claim that $V^i_k\cap V^i_{k'}=\emptyset$ whenever $V_k\cap V_{k'}=\emptyset $. If not, let there exist an $i$ and $k, k'$ such that $V^i_k\cap V^i_{k'}\neq \emptyset$. Remark \ref{translationofad}(2) gives that  $\cap_{k\in \mathbb{Z}}V^i+k\pi =V^i$ concluding that if $z\in V^i$ then $z+\pi\in V^i$. Clearly, $f(z)$ and $\pi+f(z)$ are contained in $V^{i+1}$ giving that $V^{i+1}\cap V^{i+1}+\pi \neq \emptyset$ and consequently, $\cap_{k\in \mathbb{Z}}V^{i+1}+k\pi =V^{i+1}$. Repeating the above argument at most $p$-times, we get that $V\cap V+\pi\neq \emptyset$. But this is not true.

\item The existence of infinitely many invariant Fatou components of $f_\lambda$ are proved in Lemma \ref{attractingdomaincomplete}, Lemma \ref{bakerdomaincomplete}, Lemma \ref{parabolic} and Remark  \ref{ Siegel}. These domains can be seen as an infinite number of copies (translation) of one domain.
For $|2+\lambda^2|<1$, $\lambda=i$, $2+\lambda^2=e^{2\pi i \alpha}$ for some rational number $\alpha$ and $2+\lambda^2=e^{2\pi i \alpha}$ for some bounded type irrational number $\alpha$,  the set of all (invariant) attracting domains, Baker domains, parabolic domains and Siegel discs of $f_\lambda$ respectively can be enumerated in a specific way. In fact, if $V$ is such a domain of $f_\lambda$ then for each $k\in \mathbb{Z}$, $V+k\pi$ is also such a domain, and $\{V+k\pi:k\in \mathbb{Z}\}$ is the set of all such domains.
	 
\end{enumerate}
\end{rem}

\section{The proofs}

\begin{lem}\label{lambda-pilambda}

 For each $m \in \mathbb{Z} $, the Fatou set of $f_\lambda$ is equal to the Fatou set of $f_{\lambda+m\pi}$.

\end{lem}

\begin{proof}
For $z\in \mathbb{C}$, $f_{\lambda+m\pi}(z)=\lambda+m\pi+z+\tan z=m\pi+f_\lambda(z)$ and $f^2_{\lambda+m\pi}(z)=\lambda+m\pi+m\pi+f_{\lambda}(z)+\tan (m\pi+f_{\lambda}(z))=2m\pi+f_\lambda(z)+\tan (f_\lambda(z))=2m\pi+f^2_\lambda(z)$. Hence

\begin{equation}\label{escaping}
f^n_{\lambda+m\pi}(z)=nm\pi+f^n_\lambda(z)
 \end{equation}
 
  giving that $\{f^n_\lambda\}_{n\geq 0}$ is normal at $z$ if and only if $\{f^n_{\lambda+m\pi}\}_{n\geq 0}$ is normal at $z$.
\end{proof}

\begin{lem}\label{wandeing}
For $|2+\lambda^2|\leq 1$, if $V$ is a periodic Fatou component such that $V\cap V+\pi =\emptyset$ then $V$ is a wandering domain of $f_{\lambda+m\pi}$ for $m \in \mathbb{Z}\setminus \{0\}$. Further, if $V$ is an attracting domain or a parabolic domain of $f_\lambda$ then $V$ is escaping wandering domain of $f_{\lambda+m\pi}$.
\end{lem}

\begin{proof}
Let $V$ be a $p$-periodic Fatou component of $f_\lambda$ and $\{V, V^1, V^2,..., V^{p-1}\}$ be the cycle containing $V$, where $V^i=f_\lambda^i(V)$. Since $V\cap V+\pi=\emptyset$, for $i=0,1,2,...,p-1$ using Remark \ref{translationofad}(3) we get  $V^i+k\pi \cap V^i+k'\pi=\emptyset$ for $k \neq k'$.
It follows from Lemma \ref{lambda-pilambda} that the Fatou set of $f_\lambda$ and the Fatou set of $f_{\lambda+m\pi}$ are equal. Though the Fatou set of $f_{\lambda+m\pi}$ and that of $f_{\lambda}$ are identical, the iterative behavior on their respective Fatou components differ. Now we discuss the nature of $V$ as a Fatou component of $f_{\lambda+m\pi}$. For $z\in V$, $f_{\lambda+m\pi}(z)=m\pi+f_\lambda(z) \in V^1+m\pi$ and hence, $f^n_{\lambda+m\pi}(z)=nm\pi+f^n_\lambda(z) \in V^n+nm\pi$

Since $V\cap V+\pi =\emptyset$, it follows from Remark \ref{translationofad}(2) and Remark \ref{translationofad}(3) that for each $n$, $f^n_{\lambda+m\pi}(z)$ is contained in different Fatou components and $\{V, V^1+m\pi, V^2+2m\pi, ..., V^{p-1}+m(p-1)\pi, V+mp\pi,....\}$ is the orbit of $V$ under $f_{\lambda+m\pi}$. Thus $V$ is a wandering domain of $f_{\lambda+m\pi}$.
\par 
If $V$ is an attracting domain or a parabolic domain of $f_{\lambda}$ corresponding to a $p$-periodic point $z_0$ then for each $z\in V$ 
\begin{equation}\label{escaping1}
f^{np}_{\lambda}(z)\rightarrow z_0 ~\mbox{and}~  f^{np}_{\lambda+m\pi}(z)=npm\pi+f^{np}_{\lambda}(z)\rightarrow \infty ~\mbox{as}~ n\rightarrow \infty. 
\end{equation}
For $i=0,1,2,...,p-1$ we have 
\begin{equation}\label{escaping2}
f_\lambda^i(f^{np}_{\lambda})(z)\rightarrow f_\lambda^i(z_0) ~\mbox{and}~  f_{\lambda+m\pi}^i(f^{np}_{\lambda+m\pi})(z)=(np+i)m\pi+f_\lambda^i(f^{np}_{\lambda})(z)\rightarrow \infty ~\mbox{as}~ n\rightarrow \infty. 
\end{equation}
Since $\{V, V^1+m\pi, V^2+2m\pi, ..., V^{p-1}+m(p-1)\pi, V+mp\pi,....\}$ is the orbit of $V$ under $f_{\lambda+m\pi}$, using Equation (\ref{escaping1}) and Equation (\ref{escaping2}) we have $V$ is escaping. 

\end{proof}
		
\begin{lem}\label{combine}
	For $|2+\lambda^2|<1,~\lambda=i,~2+\lambda^2=e^{2\pi i\alpha}$ for $\alpha \in \mathbb{Q}$ and $2+\lambda^2=e^{2\pi i\alpha}$ for some bounded type irrational $\alpha$, the Fatou set of $f_{\lambda+m\pi}$ has $m$ many escaping wandering domains with disjoint grand orbit.
\end{lem}

\begin{proof}
	It is evident from Lemma \ref{lambda-pilambda} that the Fatou set of $f_\lambda$ and $f_{\lambda+m\pi}$ are identical for $m\in \mathbb{Z}$. In Spite of having identical Fatou sets, the internal behavior of the Fatou components of $f_\lambda$ and that of $f_{\lambda+m\pi}$ are poles apart.
	Lemma \ref{cifcupperhalf} indicates the existence of a completely invariant Baker domain of $f_\lambda$ and $f_{\lambda+m\pi}$ that contains the upper-half plane. As both of them contain the upper-half plane and are completely invariant Baker domain then they are the same Fatou component having a similar property. For $|2+\lambda^2|<1,~\lambda=i,~2+\lambda^2=e^{2\pi i\alpha}$ for $\alpha \in \mathbb{Q}$ and $2+\lambda^2=e^{2\pi i\alpha}$ for some bounded type irrational $\alpha$, it follows from  Lemma~\ref{attractingdomaincomplete}, Lemma \ref{bakerdomaincomplete}, Lemma \ref{parabolic} and Remark \ref{ Siegel} that along with the primary Fatou component, the Fatou set of  $f_{\lambda}$ contains infinitely many invariant attracting domains, Baker domains, parabolic domains and  Siegel discs respectively. For a given such attracting domain, Baker domain, parabolic domain, and  Siegel disc $V$, all other attracting domains, Baker domains, parabolic domains, and  Siegel discs respectively are of the form $V+k\pi$, $k\in \mathbb{Z}$ (Remark \ref{translationofad}(4)). Since each such Fatou components are disjoint from each other, $V+k\pi\cap V+k'\pi=\emptyset$ for all $k \neq k'$. Lemma \ref{wandeing} concludes that $V$ is a wandering Fatou component of $f_{\lambda+m\pi}$. As $f_{\lambda+m\pi}$ maps $V$ into $V+m\pi$ (Remark \ref{translationofad}), $V+\pi$ is mapped into $V+(m+1)\pi$. Hence $f_{\lambda+m\pi}$ has $|m|$ many wandering domains with different grand orbits, namely $V, V+\pi, V+2 \pi,\cdots, V+(m-1)\pi$. It now follows from Equation (\ref{escaping1}) that $V$ is escaping.
	
\end{proof}

Now we are ready to prove Theorem \ref{wanderingattracting}.

\begin{proof}
	Let $V$ be an invariant attracting domain of $f_\lambda$. Using Lemma \ref{combine}  we get that $f_{\lambda+m\pi}$ has $|m|$ many wandering domains with different grand orbits.
	
	\begin{enumerate}
		\item It follows from the second part of Lemma \ref{wandeing} that these wandering domains are escaping.
		\item Note that  $V+n\pi$ for  each $n \in \mathbb{Z}$, contains a critical point of $f_{\lambda}$ and the set of all the critical points of $f_{\lambda}$ and that of $f_{\lambda+m\pi}$ are identical. Thus each such wandering domain $V+n\pi$ contains a critical point of $f_{\lambda+m\pi}$. Since each of the attracting domains of $f_\lambda$ contains exactly one critical point, each such wandering domain of $f_{\lambda+m\pi}$ contains exactly one critical point. By Lemma~\ref{attractingdomaincomplete}, $f_{\lambda+m\pi}: V+n \pi \to V + (n+m) \pi$ is proper. Its degree is $2$ by the Riemann-Hurwitz formula. Taking $ W_n =V+n \pi$ and $W_{n+1}=V+ (n+m) \pi$, it is seen that $f_{\lambda+m\pi}: W_n \to W_{n+1}$ is a proper map with degree $2$.
		\item  If $W'$ is a wandering domain in the grand orbit of $W(=W_0=V)$ and is different from every $W_n$ then there is no critical point in $W'$ and the map $f_{\lambda+m\pi}$ is one-one on $W'$ by the Riemann-Hurwitz formula.
		\item Let $z_0\in V$ be the attracting fixed point of $f_\lambda$. Clearly, $z_0+n\pi$ is also an attracting fixed point of $f_\lambda$ contained in $V+n\pi$ for each $n \in \mathbb{Z}$. Since $V$ is an attracting domain, for each $z\in V $ and $\epsilon>0$, there exists an $n_o$ such that $\overline{D_\epsilon(z_0)}\subset V$ and $f^n_\lambda(z) \in D_\epsilon(z_0)$ for all $n\geq n_0$. Here $D_\epsilon(z_0)$ is the ball of radius $\epsilon$ centered at $z_0$. Thus $\liminf_{n\rightarrow \infty}dist(f_\lambda^n(z),\partial V)>0$ for all $z\in V$.  Note that $\overline{D_\epsilon(z_0+n\pi)}\subset V+n\pi$. It follows from Equation~(\ref{escaping}) that for a given $\epsilon>0$ and $z \in V$ if  $f_{\lambda}^n (z) \to z_0$ as $n \to \infty$ then $f_{\lambda+m\pi}^n (z) \to \infty$ and $|f_{\lambda+m\pi}^n (z)-\pi n m- z_0|<\epsilon$ as $n \to \infty$. This implies that $\liminf_{n\rightarrow \infty}dist(f_{\lambda+m\pi}^n(z),\partial (V+nm\pi))>0$ for all $z\in V$. Again taking $ W =V$ and $W_{n}=V+nm \pi$ we have $\liminf_{n\rightarrow \infty}dist(f_\lambda^n(z),\partial W_n)>0$ for all $z\in W$.
	\end{enumerate}
	It follows from Lemma \ref{attractingdomaincomplete} that $f_{\lambda+m\pi}$ does not have any other periodic Fatou component or any other wandering domain. This completes the proof.
\end{proof}

We now authenticate Theorem \ref{wanderingbaker}.

\begin{proof}
	Let $V$ be an invariant Baker domain of $f_i$. Using Lemma \ref{combine}  we get that $f_{i+m\pi}$ has $|m|$ many wandering domains with different grand orbits.
	
	\begin{enumerate}
		\item It also follows from Lemma \ref{combine} that these wandering domains are escaping.
		\item Note that  $V+n\pi$ for  each $n \in \mathbb{Z}$, contains a critical point of $f_{i}$ and the set of all the critical points of $f_{i}$ and that of $f_{i+m\pi}$ are identical. Thus each such wandering domain $V+n\pi$ contains a critical point of $f_{i+m\pi}$. Since each of the Baker domains of $f_i$ contains exactly one critical point, each such wandering domain of $f_{i+m\pi}$ contains exactly one critical point. By Lemma~\ref{attractingdomaincomplete}, $f_{i+m\pi}: V+n \pi \to V + (n+m) \pi$ is proper. Its degree is $2$ by the Riemann-Hurwitz formula. Taking $ W_n =V+n \pi$ and $W_{n+1}=V+ (n+m) \pi$, it is seen that $f_{i+m\pi}: W_n \to W_{n+1}$ is a proper map with degree $2$.
		\item  If $W'$ is a wandering domain in the grand orbit of $W(=W_0=V)$ and is different from every $W_n$ then there is no critical point in $W'$ and the map $f_{i+m\pi}$ is one-one on $W'$ by the Riemann-Hurwitz formula.
		\item Note that $\{\Re(z): z \in V\}$ is bounded. Since for each $z\in V$, $f_i^n(z)\rightarrow \infty$ as $n\rightarrow \infty$, $\Im(f_i^n(z))\rightarrow -\infty$ whenever $n \rightarrow \infty$. Thus using Equation~(\ref{escaping}) we have $\{\Re(f^n_{i+m\pi}(z))\rightarrow \infty\}$ and $\{\Im(f^n_{i+m\pi}(z))\rightarrow -\infty\}$. By considering $W=V$, this gives that
		$\liminf_{n\rightarrow \infty}dist(f_{i+m\pi}^n(z),\partial W_n)=0$ and, $\Im(f^n_{i+m\pi}(z))\rightarrow -\infty$ for all $z\in W$ whenever $n \rightarrow \infty$.
		
	\end{enumerate}
	It follows from Lemma \ref{bakerdomaincomplete} that $f_{\lambda+m\pi}$ does not have any other periodic Fatou component or any other wandering domain. This completes the proof.
\end{proof}

The proof Theorem \ref{wanderingparabolic} is portrayed here.

\begin{proof}
	Let $V$ be an invariant parabolic domain of $f_\lambda$. Using Lemma \ref{combine}  we get that $f_{\lambda+m\pi}$ has $|m|$ many wandering domains with different grand orbits.
	
	\begin{enumerate}
		\item It also follows from Lemma \ref{combine} that these wandering domains are escaping.
		
		\item Using the similar argument as discussed in Theorem \ref{wanderingattracting}(2), it can be proved that  $f_{i+m\pi}: W_n \to W_{n+1}$ is a proper map with degree $2$.
		
		\item  If $W'$ is a wandering domain in the grand orbit of $W(=W_0=V)$ and is different from every $W_n$ then there is no critical point in $W'$ and the map $f_{\lambda+m\pi}$ is one-one on $W'$ by the Riemann-Hurwitz formula.
		
		\item Let $z_0\in \partial V$ be the parabolic fixed point of $f_\lambda$. Clearly, $z_0+n\pi$ is also a parabolic fixed point of $f_\lambda$ contained in $\partial (V+n\pi)$ for each $n \in \mathbb{Z}$. Since $V$ is a parabolic domain, for each $z\in V $ and $\epsilon>0$, there exists an $n_o$ such that  $f^n_\lambda(z) \in B_\epsilon(z_0)$ for all $n\geq n_0$ and hence, $\liminf_{n\rightarrow \infty}dist(f_\lambda^n(z),\partial V)=0$. It follows from Equation~(\ref{escaping}) that for a given $\epsilon>0$ and $z \in V$, $f_{\lambda+m\pi}^n (z) \to \infty$ and $|f_{\lambda+m\pi}^n (z)-\pi n m- z_0|<\epsilon$ as $n \to \infty$. This implies that $\liminf_{n\rightarrow \infty}dist(f_{\lambda+m\pi}^n(z),\partial (V+nm\pi))=0$ for all $z\in V$. Again taking $ W =V$ and $W_{n}=V+nm \pi$ we have $\liminf_{n\rightarrow \infty}dist(f_\lambda^n(z),\partial W_n)=0$ for all $z\in W$.
	\end{enumerate}
	It follows from Lemma \ref{parabolic} that $f_{\lambda+m\pi}$ does not have any other attracting domain, parabolic domain,  Siegel disc, or Herman ring.
\end{proof}

Following is the proof of Theorem \ref{wanderingsiegel}.

\begin{proof}
	Let $V$ be an invariant  Siegel disc of $f_\lambda$. Using Lemma \ref{combine}  we get that $f_{\lambda+m\pi}$ has $|m|$ many wandering domains with different grand orbits.
	
	\begin{enumerate}
		\item It also follows from Lemma \ref{combine} that these wandering domains are escaping.
		
		\item Since $f_\lambda$ has a  Siegel disc, $f_\lambda$ has at least one critical point in the Julia set. It follows from Lemma \ref{criticalpoint}(1) and Remark \ref{translationofad}(1) that each critical point of $f_\lambda$ is in the Julia set. Note that the critical points of $f_\lambda$ and that of $f_{\lambda+m\pi}$ are the same. It now follows from Lemma \ref{lambda-pilambda} that all of the critical points of $f_{\lambda+m\pi}$ are in the Julia set. Using the Riemann-Hurwitz formula we get   $f_{\lambda+m\pi}: W \to W_{1}$ is a proper map with degree $1$.

		\item Since $V$ is a  Siegel disc, for each $z\in V$, there exists an $f_\lambda$-invariant Jordan curve $\gamma$ containing $z$ such that $dist(\gamma,\partial V)>0$ and hence, $dist(f_\lambda^n(z),\partial V)>0$. In particular, $\lim\inf_{n\rightarrow \infty} dist(f^n_\lambda(z), \partial V)>0$. Again by the definition of  Siegel disc, we have $dist(f^n_{\lambda}(z),f^n_{\lambda}(z'))=dist(z,z')$ for all $z,z'\in V$. Taking $ W =V$ and $W_{n}=V+nm \pi$ and using  Equation (\ref{escaping}) we get  $\lim\inf_{n\rightarrow \infty} dist(f^n_{\lambda+m\pi}(z), \partial W_n)>0$ and $dist(f^n_{\lambda+m\pi}(z),f^n_{\lambda+m\pi}(z'))=dist(z,z')$ for all $z,z'\in W$ .
	\end{enumerate}
	The existence of a completely invariant Baker domain of $f_{\lambda+m\pi}$ concludes the non-existence of the Herman ring of $f_{\lambda+m\pi}$. Since all the singular values of $f_{\lambda+m\pi}$ are contained in the Julia set,  $f_{\lambda+m\pi}$ does not have any other attracting domain or parabolic domain.
\end{proof}
\newpage

\bibliographystyle{amsplain}
\addcontentsline{toc}{chapter}{\numberline{}References}
 
\end{document}